\documentstyle{article}

\def\prop#1#2{\vspace{2ex} \noindent\textsc{#1.} \emph{#2} \par \vspace{2ex}}
\def\mn{\!-\!}

\begin{document}

\title{\textbf{A Planarity Criterion for Graphs}}
\author{\textsc{Kosta Do\v sen} and \textsc{Zoran Petri\' c}
\\[1ex]
{\small Mathematical Institute, SANU}\\[-.5ex]
{\small Knez Mihailova 36, p.f.\ 367, 11001 Belgrade,
Serbia}\\[-.5ex]
{\small email: \{kosta, zpetric\}@mi.sanu.ac.rs}}
\date{}
\maketitle

\vspace{-4ex}

\begin{abstract}
\noindent It is proven that a connected graph is planar if and
only if all its cocycles with at least four edges are ``grounded''
in the graph. The notion of grounding of this planarity criterion,
which is purely combinatorial, stems from the intuitive idea that
with planarity there should be a linear ordering of the edges of a
cocycle such that in the two subgraphs remaining after the removal
of these edges there can be no crossing of disjoint paths that
join the vertices of these edges. The proof given in the paper of
the right-to-left direction of the equivalence is based on
Kuratowski's Theorem for planarity involving $K_{3,3}$ and $K_5$,
but the criterion itself does not mention $K_{3,3}$ and $K_5$.
Some other variants of the criterion are also shown necessary and
sufficient for planarity.
\end{abstract}


\section{Introduction}
In this note we prove the necessity and sufficiency of a rather
simple planarity criterion for graphs (which, as usual, and as in
\cite{H69}, are understood to be finite, not directed, without
multiple edges and without loops). We prove that a connected graph
is planar if and only if all its cocycles with at least four edges
are \emph{grounded} in the graph. The notion of grounding, which
is purely combinatorial, will be defined precisely later in this
section. Our planarity criterion, which is based on this notion,
stems from the intuitive idea that with planarity there should be
a linear ordering of the edges of a cocycle such that in the two
subgraphs remaining after the removal of these edges there can be
no crossing of disjoint paths that join the vertices of these
edges. The criterion will become clear with the examples of the
next section. As far as we know, this criterion is new, and it is
formulated without mentioning the graphs $K_{3,3}$ and $K_5$ of
Kuratowski's planarity criterion (see \cite{K30} and \cite{H69},
Chapter 11).

The proof of necessity for this criterion, i.e.\ of the
left-to-right direction of the equivalence, is easy. The proof
given here of sufficiency, i.e.\ of the remaining direction,
relies however on Kuratowski's Theorem on planarity. We use
actually the sufficiency direction for Kuratowski's criterion,
which is more difficult to prove than the necessity direction. We
suppose that an independent proof could be given for the
sufficiency of our criterion, but we do not expect it to be
shorter than the proof of sufficiency for Kuratowski's criterion,
and, as far as we can see, it would rely on similar ideas. So we
do not find it worthwhile to go into such a new proof, which would
make the paper longer.

For newer papers giving, like ours, planarity criteria for graphs
alternative to Kuratowski's one may consult the references of
\cite{LS10} and \cite{M11}. An extensive bibliography for various
matters related to Kuratowski's Theorem may be found
in~\cite{T81}.

We will follow the terminology and notation of \cite{H69} whenever
we can, except that instead of \emph{point} and \emph{line} we use
respectively \emph{vertex} and \emph{edge} for abstract graphs
too, and not only for geometric graphs (embedded in
$\mathbf{R}^3$)---this usage is presumably more common.

A \emph{cutset} of a connected graph is a set of its edges whose
removal (see \cite{H69}, Chapter~2) results in a disconnected
graph, and a \emph{cocycle} is defined in \cite{H69} (Chapter~4)
as a minimal cutset (note that elsewhere, as e.g.\ in \cite{LS10},
the terminology might be different). We call a cocycle \emph{big}
when it has at least four edges. Given a cocycle $C$ of a
connected graph $G$, let $G'$ and $G''$ be the two connected
subgraphs of $G$ obtained by removing the edges of $C$ from $G$.
We keep this convention throughout the paper.

Consider four distinct edges $x_1$, $x_2$, $x_3$ and $x_4$ of a
big cocycle $C$ of $G$; we assume that for $i\in\{1,\dots,4\}$ we
have that $x_i$ is $u_iv_i$, and the vertices $u_1$, $u_2$, $u_3$
and $u_4$ are in $G'$, while the vertices $v_1$, $v_2$, $v_3$ and
$v_4$ are in $G''$. Note that although the four edges $x_1$,
$x_2$, $x_3$ and $x_4$ are distinct, the vertices $u_1$, $u_2$,
$u_3$ and $u_4$ need not be distinct, and the same for the
vertices $v_1$, $v_2$, $v_3$ and~$v_4$.

We say that $\{x_1,x_3\}$ and $\{x_2,x_4\}$ are \emph{disparate}
in $G'$ when in $G'$ we have a ${u_1\mn u_3}$ path and a ${u_2\mn
u_4}$ path with no vertex in common. Analogously, we say that
$\{x_1,x_3\}$ and $\{x_2,x_4\}$ are \emph{disparate} in $G''$ when
in $G''$ we have a ${v_1\mn v_3}$ path and a ${v_2\mn v_4}$ path
with no vertex in common. (Something analogous to our notion of
disparate pairs of edges is given for vertices in the notion of
skew $C$-components; see \cite{T81}, Section~2.)

In a sequence $X_1a_1\dots X_na_nX_{n+1}$, where $n\geq 1$, the
sequence $a_1\ldots a_n$ is a nonempty \emph{subsequence}; here
$a_i$, for $i\in\{1,\dots,n\}$, is a member of our sequence and
$X_j$, for $j\in\{1,\dots,n\! +\! 1\}$, is a sequence, possibly
empty, of such members.

A big cocycle $C$ of $G$ is \emph{grounded} in $G$ when there is a
sequence without repetitions containing all its edges such that for every
subsequence $x_1x_2x_3x_4$ of this sequence we have that
$\{x_1,x_3\}$ and $\{x_2,x_4\}$ are disparate neither in $G'$ nor
in~$G''$.

The theorem giving our planarity criterion is the following.

\prop{Theorem}{A connected graph is planar if and only if each of
its big cocycles is grounded in it.}

An arbitrary graph is planar when, of course, each of its
connected subgraphs is planar. So this theorem yields easily a
planarity criterion for arbitrary graphs.

In the next two sections we consider preliminary matters, which we
use in Section~4 to give a proof of our theorem. At the end of
that section, and at the very end of the paper, we envisage some
variants of our theorem, which are easily derived from our proof,
and which may be interesting from an algorithmic point of view.

\section{The big cocycles of $K_{3,3}$ and $K_5$}

It happens that if $G$ is $K_{3,3}$ or $K_5$, then all the big
cocycles of $G$ are not grounded. For the proof of the Theorem it
is however enough that at least one of these big cocycles of $G$
is not grounded.

If $G$ is $K_{3,3}$ we have just two types of big cocycles. The
first type is given by the dotted edges in the following picture
of~$G$:
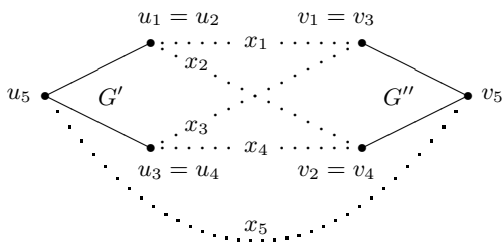
\begin{center}
\begin{picture}(160,100)(0,-10)

\put(0,50){\circle*{3}} \put(40,30){\circle*{3}}
\put(40,70){\circle*{3}} \put(120,30){\circle*{3}}
\put(120,70){\circle*{3}} \put(160,50){\circle*{3}}

\put(0,50){\line(2,-1){40}} \put(0,50){\line(2,1){40}}
\put(160,50){\line(-2,-1){40}} \put(160,50){\line(-2,1){40}}

{\thicklines \qbezier[35](0,50)(80,-60)(160,50)}

\put(80,70){\small\makebox(0,0){$x_1$}}
\put(57.5,62){\small\makebox(0,0){$x_2$}}
\put(57.5,38){\small\makebox(0,0){$x_3$}}
\put(80,30){\small\makebox(0,0){$x_4$}}
\put(80,0){\small\makebox(0,0){$x_5$}}

\put(35,79){\small\makebox(0,0)[l]{$u_1=u_2$}}
\put(125,79){\small\makebox(0,0)[r]{$v_1=v_3$}}
\put(35,21){\small\makebox(0,0)[l]{$u_3=u_4$}}
\put(125,21){\small\makebox(0,0)[r]{$v_2=v_4$}}
\put(-5,50){\small\makebox(0,0)[r]{$u_5$}}
\put(165,50){\small\makebox(0,0)[l]{$v_5$}}

\multiput(40,30)(5,0){7}{\makebox(0,0){\circle*{.3}}}
\multiput(90,30)(5,0){6}{\makebox(0,0){\circle*{.3}}}

\multiput(40,70)(5,0){7}{\makebox(0,0){\circle*{.3}}}
\multiput(90,70)(5,0){6}{\makebox(0,0){\circle*{.3}}}

\multiput(40,30)(5,2.5){3}{\makebox(0,0){\circle*{.3}}}
\multiput(65,42.5)(5,2.5){11}{\makebox(0,0){\circle*{.3}}}

\multiput(40,70)(5,-2.5){3}{\makebox(0,0){\circle*{.3}}}
\multiput(65,57.5)(5,-2.5){11}{\makebox(0,0){\circle*{.3}}}

\put(20,50){\small\makebox(0,0)[l]{$G'$}}
\put(140,50){\small\makebox(0,0)[r]{$G''$}}

\end{picture}
\end{center}
For $i\in\{1,\ldots,5\}$, we have that $x_i$ is the edge $u_iv_i$,
and analogously in the other pictures below.

This cocycle of $G$ is not grounded in $G$. For example, if we
take the sequence $x_1x_2x_3x_4x_5$ of the edges of our big
cocycle, then for the subsequence $x_1x_2x_3x_4$ we have that
$\{x_1,x_3\}$ and $\{x_2,x_4\}$ are disparate in $G''$, since the
one-vertex paths $v_1$ and $v_2$ have no vertex in common.

As another example, take the sequence $x_1x_2x_5x_4x_3$. Then for
the subsequence $x_1x_2x_5x_4$ we have that $\{x_1,x_5\}$ and
$\{x_2,x_4\}$ are disparate in $G''$, since the path $v_1v_5$ and
the one-vertex path $v_2$ have no vertex in common.

Up to renaming of indices, the sequences in these two examples are
the only two different sorts of sequences with our first type of
cocycle for $K_{3,3}$. In our cocycle we have two kinds of edges:
on the one hand, $x_1$, $x_2$, $x_3$ and $x_4$, each adjacent on
both ends to another edge of the cocycle, and on the other hand
$x_5$, adjacent on both ends to no other edge of the cocycle. In
the first example, in $\{x_1,x_3\}$ and $\{x_2,x_4\}$ we have only
edges of the first kind. In the second example, in $\{x_1,x_5\}$
and $\{x_2,x_4\}$ we have also the edge $x_5$ of the second kind.
Since $G'$ and $G''$ are isomorphic graphs, this exhausts all
possibilities.

The second type of big cocycle with $G$ being $K_{3,3}$ is given
by the dotted edges in the following picture of~$G$:
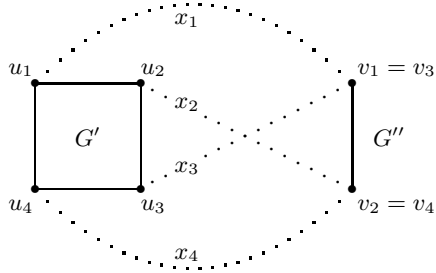
\begin{center}
\begin{picture}(120,110)(0,-5)

\put(0,30){\circle*{3}} \put(0,70){\circle*{3}}
\put(40,30){\circle*{3}} \put(40,70){\circle*{3}}
\put(120,30){\circle*{3}} \put(120,70){\circle*{3}}

\put(0,30){\line(1,0){40}} \put(0,30){\line(0,1){40}}
\put(40,70){\line(-1,0){40}} \put(40,70){\line(0,-1){40}}
\put(120,70){\line(0,-1){40}}

{\thicklines \qbezier[25](0,30)(60,-30)(120,30)
\qbezier[25](0,70)(60,130)(120,70)}

\put(57.5,94){\small\makebox(0,0){$x_1$}}
\put(57.5,62){\small\makebox(0,0){$x_2$}}
\put(57.5,38){\small\makebox(0,0){$x_3$}}
\put(57.5,5){\small\makebox(0,0){$x_4$}}

\put(-10,75){\small\makebox(0,0)[l]{$u_1$}}
\put(40,75){\small\makebox(0,0)[l]{$u_2$}}
\put(122,75){\small\makebox(0,0)[l]{$v_1=v_3$}}
\put(40,23.5){\small\makebox(0,0)[l]{$u_3$}}
\put(-10,23.5){\small\makebox(0,0)[l]{$u_4$}}
\put(122,23.5){\small\makebox(0,0)[l]{$v_2=v_4$}}

\multiput(40,30)(5,2.5){3}{\makebox(0,0){\circle*{.3}}}
\multiput(65,42.5)(5,2.5){11}{\makebox(0,0){\circle*{.3}}}

\multiput(40,70)(5,-2.5){3}{\makebox(0,0){\circle*{.3}}}
\multiput(65,57.5)(5,-2.5){11}{\makebox(0,0){\circle*{.3}}}

\put(20,50){\small\makebox(0,0){$G'$}}
\put(140,50){\small\makebox(0,0)[r]{$G''$}}

\end{picture}
\end{center}

This cocycle of $G$ is not grounded in $G$. For example, if we
take the sequence $x_1x_2x_3x_4$ of the edges of our big cocycle,
then for the subsequence $x_1x_2x_3x_4$ (which is our sequence
itself) we have that $\{x_1,x_3\}$ and $\{x_2,x_4\}$ are disparate
in $G''$, since the one-vertex paths $v_1$ and $v_2$ have no
vertex in common.

As another example, take the sequence $x_1x_3x_2x_4$. Then for the
subsequence $x_1x_3x_2x_4$ we have that $\{x_1,x_2\}$ and
$\{x_3,x_4\}$ are disparate in $G'$, since the paths $u_1u_2$ and
$u_3u_4$ have no vertex in common.

Since the edges in our cocycle are of the same kind (unlike what
we had with the previous type of cocycle, with five edges), up to
renaming of indices the sequences in these two examples are the
only kind of sequences with our second type of cocycle for
$K_{3,3}$. The cocycles that have $G'$ with five vertices and
$G''$ with a single vertex have three edges, and are hence not
big. There are no other types of cocycle for~$K_{3,3}$.

If $G$ is $K_5$, then we have just two types of cocycles, and they
are both big. The first type is given by the dotted edges in the
following picture of~$G$:
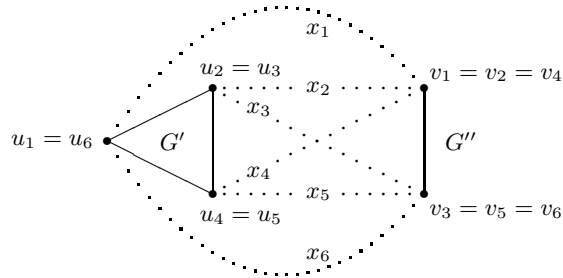
\begin{center}
\begin{picture}(120,115)(0,-5)

\put(0,50){\circle*{3}} \put(40,30){\circle*{3}}
\put(40,70){\circle*{3}} \put(120,30){\circle*{3}}
\put(120,70){\circle*{3}}

\put(0,50){\line(2,-1){40}} \put(0,50){\line(2,1){40}}
\put(40,30){\line(0,1){40}} \put(120,30){\line(0,1){40}}

{\thicklines \qbezier[30](0,50)(60,-40)(120,30)
\qbezier[30](0,50)(60,140)(120,70)}

\put(80,92){\small\makebox(0,0){$x_1$}}
\put(80,70){\small\makebox(0,0){$x_2$}}
\put(57.5,62){\small\makebox(0,0){$x_3$}}
\put(57.5,38){\small\makebox(0,0){$x_4$}}
\put(80,30){\small\makebox(0,0){$x_5$}}
\put(80,7){\small\makebox(0,0){$x_6$}}

\put(35,77){\small\makebox(0,0)[l]{$u_2=u_3$}}
\put(122,75){\small\makebox(0,0)[l]{$v_1=v_2=v_4$}}
\put(35,21){\small\makebox(0,0)[l]{$u_4=u_5$}}
\put(122,23.5){\small\makebox(0,0)[l]{$v_3=v_5=v_6$}}
\put(-5,50){\small\makebox(0,0)[r]{$u_1=u_6$}}

\multiput(40,30)(5,0){7}{\makebox(0,0){\circle*{.3}}}
\multiput(90,30)(5,0){6}{\makebox(0,0){\circle*{.3}}}

\multiput(40,70)(5,0){7}{\makebox(0,0){\circle*{.3}}}
\multiput(90,70)(5,0){6}{\makebox(0,0){\circle*{.3}}}

\multiput(40,30)(5,2.5){3}{\makebox(0,0){\circle*{.3}}}
\multiput(65,42.5)(5,2.5){11}{\makebox(0,0){\circle*{.3}}}

\multiput(40,70)(5,-2.5){3}{\makebox(0,0){\circle*{.3}}}
\multiput(65,57.5)(5,-2.5){11}{\makebox(0,0){\circle*{.3}}}

\put(20,50){\small\makebox(0,0)[l]{$G'$}}
\put(140,50){\small\makebox(0,0)[r]{$G''$}}

\end{picture}
\end{center}

This cocycle of $G$ is not grounded in $G$. For example, if we
take the sequence $x_1x_2x_3x_4x_5x_6$ of the edges of our big
cocycle, then for the subsequence $x_2x_3x_4x_5$ we have that
$\{x_2,x_4\}$ and $\{x_3,x_5\}$ are disparate in $G''$, since the
one-vertex paths $v_1$ and $v_3$ have no vertex in common.

As another example, take the sequence $x_1x_2x_4x_3x_5x_6$. Then
for the subsequence $x_2x_4x_3x_5$ we have that $\{x_2,x_3\}$ and
$\{x_4,x_5\}$ are disparate in $G'$, since the one-vertex paths
$u_2$ and $u_4$ have no vertex in common. Since all the edges in
our cocycle are of the same kind, up to renaming of indices the
sequences in these two examples are the only kinds of sequences
with this first type of cocycle for~$K_5$.

The second type of big cocycle with $G$ being $K_5$ is given by
the dotted edges in the following picture of~$G$:
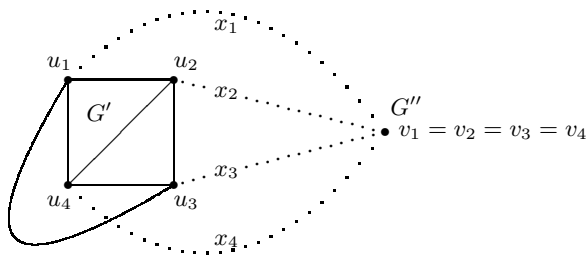
\begin{center}
\begin{picture}(120,96)
\put(0,30){\circle*{3}} \put(0,70){\circle*{3}}
\put(40,30){\circle*{3}} \put(40,70){\circle*{3}}
\put(120,50){\circle*{3}}

\put(0,30){\line(1,0){40}} \put(0,30){\line(0,1){40}}
\put(40,70){\line(-1,0){40}} \put(40,70){\line(0,-1){40}}
\put(0,30){\line(1,1){40}}

\qbezier(0,70)(-60,-30)(40,30) {\thicklines
\qbezier[25](0,30)(60,-30)(120,50)
\qbezier[25](0,70)(60,130)(120,50)}

\put(60,90){\small\makebox(0,0){$x_1$}}
\put(60,65){\small\makebox(0,0){$x_2$}}
\put(60,35){\small\makebox(0,0){$x_3$}}
\put(60,9){\small\makebox(0,0){$x_4$}}

\put(-8,76){\small\makebox(0,0)[l]{$u_1$}}
\put(40,76){\small\makebox(0,0)[l]{$u_2$}}
\put(40,23.5){\small\makebox(0,0)[l]{$u_3$}}
\put(-8,23.5){\small\makebox(0,0)[l]{$u_4$}}
\put(125,50){\small\makebox(0,0)[l]{$v_1=v_2=v_3=v_4$}}

\multiput(40,30)(4,1){4}{\makebox(0,0){\circle*{.3}}}
\multiput(68,37)(4,1){13}{\makebox(0,0){\circle*{.3}}}

\multiput(40,70)(4,-1){4}{\makebox(0,0){\circle*{.3}}}
\multiput(68,63)(4,-1){13}{\makebox(0,0){\circle*{.3}}}

\put(12,58){\small\makebox(0,0){$G'$}}
\put(122,60){\small\makebox(0,0)[l]{$G''$}}

\end{picture}
\end{center}

This cocycle of $G$ is not grounded in $G$. For example, if we
take the sequence $x_1x_2x_3x_4$ of the edges of our big cocycle,
then for the subsequence $x_1x_2x_3x_4$ we have that $\{x_1,x_3\}$
and $\{x_2,x_4\}$ are disparate in $G'$, since the paths $u_1u_3$
and $u_2u_4$ have no vertex in common. In $G'$ we have also that
$\{x_1,x_2\}$ and $\{x_3,x_4\}$ are disparate, because the paths
$u_1u_2$ and $u_3u_4$ have no vertex in common, while
$\{x_1,x_4\}$ and $\{x_2,x_3\}$ are disparate because the paths
$u_1u_4$ and $u_2u_3$ have no vertex in common. Since all the
edges in our cocycle are of the same kind, up to renaming of
indices the sequence in our example is the only kind of sequence
with our second type of cocycle for $K_5$. There are no other
types of cocycle for~$K_5$.

\section{Extending cocycles}
Let $G$ be a connected graph, let $H$ be a connected subgraph of
$G$, and let $D$ be a cocycle of $H$. We will prove the following.

\prop{Lemma}{There is a cocycle $C$ of $G$ such that $D\subseteq
C$. If $D$ is a big cocycle not grounded in $H$, then $C$ is a big
cocycle not grounded in~$G$.}

\noindent\textsc{Proof.} Let $J$ be the induced subgraph (see
\cite{H69}, Chapter~2, for this notion) of $G$ with the same set
of vertices as $H$. Since $H$ is connected, $J$ must be connected
too. We define first a cocycle $E$ of $J$ such that $D\subseteq
E$. The cocycle $E$ will be the set of all edges $uv$ of $J$ such
that $u$ is in $H'$ and $v$ is in $H''$, where $H'$ and $H''$ are
the connected subgraphs of $H$ obtained by removing $D$ from~$H$.

The remainder of the proof will be made by induction on the number
$n$ of vertices in $G$ that are not in $H$. If $n=0$, then $G$ and
$H$ have the same sets of vertices, and $G$ and $J$ coincide. The
cocycle $C$ of the lemma will then be~$E$.

Suppose for the induction hypothesis that $K$ is an induced
subgraph of $G$, that $K$ is connected, that $F$ is a cocycle of
$K$ such that $D\subseteq F$, and let there be a vertex of $G$ not
in $K$. Let the removal of $F$ from $K$ result in the connected
subgraphs $K'$ and $K''$ of~$K$.

For every vertex $u$ in $G$ that is not in $K$ consider the set
$L_u'$ of edges $uv$ of $G$ with $v$ a vertex of $K'$; the set
$L_u''$ is defined in the same manner with respect to $K''$. There
must be a vertex $u$ such that $L_u'\cup L_u''\neq\emptyset$,
because $G$ is connected.

Let $M$ be the graph obtained by adding to $K$ such a vertex $u$
and all the edges in $L_u'\cup L_u''$. We obtain a cocycle $N$ of
$M$ by stipulating that $N$ is $F\cup L_u'$ or $F\cup L_u''$ if
$L_u'\neq\emptyset$ and $L_u''\neq\emptyset$; otherwise (i.e., if
$L_u'=\emptyset$ or $L_u''=\emptyset$), we have that $N$ is~$F$.

It is clear that $M$ is an induced subgraph of $G$, that $M$ is
connected, and that $N$ is a cocycle of $M$ such that $D\subseteq
N$. So by induction we conclude that there is a cocycle $C$ of $G$
such that $D\subseteq C$.

If $D$ is big, then $C$ is, of course, also big. If $D$ is not
grounded in $H$, then for an arbitrary sequence of its members we
have a subsequence $x_1x_2x_3x_4$ such that $\{x_1,x_3\}$ and
$\{x_2,x_4\}$ are disparate either in $H'$ or in $H''$. It is easy
to conclude that $\{x_1,x_3\}$ and $\{x_2,x_4\}$ are hence
disparate either in $G'$ or in $G''$, which are the subgraphs of
$G$ obtained by removing $C$ from $G$. From that we conclude
easily that $C$ is not grounded in $G$. This concludes the proof
of the Lemma.

\section{Proof of the Theorem}
For the proof of the right-to-left direction, i.e.\ the
sufficiency direction, suppose a connected graph $G$ is not
planar. By Kuratowski's Theorem, there must be a subgraph $H$ of
$G$ that is homeomorphic to either $K_{3,3}$ or $K_5$. Since edge subdivision
produces out of a cocycle that is not grounded
in a graph at least one cocycle that is not grounded in the graph that
results from the subdivision, we have that what is shown in
Section~2 holds also for every graph $H$ homeomorphic to $K_{3,3}$ or $K_5$.
Hence there is a big cocycle $D$ of $H$
that is not grounded in $H$. By the Lemma of Section~3, there is a
big cocycle $C$ of $G$ that is not grounded in~$G$.

For the proof of the left-to-right direction, i.e.\ the necessity
direction, suppose $G$ is a planar graph, with $\Gamma$ being a
plane graph realizing $G$. Let $\Gamma^\ast$ be the geometric dual
(see \cite{H69}, Chapter 11, for this notion) of $\Gamma$. For
every big cocycle $C$ of $G$, there is a cocycle $K$ of $\Gamma$
realizing $C$, and a cycle $K^\ast$ of $\Gamma^\ast$ such that the
edges of $K$ and the edges of $K^\ast$ correspond bijectively to
each other by intersecting in a single point (see \cite{B73},
Theorem~3, Section 2.2). The cycle $K^\ast$ gives a sequence of
the members of $K$, and hence a sequence of the members of~$C$.

For every subsequence $x_1x_2x_3x_4$ of this sequence we have that
$\{x_1,x_3\}$ and $\{x_2,x_4\}$ are disparate neither in $G'$ nor
in $G''$. Otherwise, we would have in $\Gamma$ two intersecting
paths without common vertex (see Lemma~2 of \cite{M84}). So $C$ is
grounded in $G$. This concludes the proof of the Theorem.

\vspace{2ex}

There are necessary and sufficient conditions for planarity that
are variants of the criterion in our theorem. One obtains these
other criteria by restricting the big cocycles mentioned in the
Theorem. For example, one may restrict them to big cocycles with
at least five edges. This is because, as shown in Section~2, for
$G$ being either $K_{3,3}$ or $K_5$ there is one such big cocycle
of $G$ not grounded in $G$. It is clear for the Lemma of Section~3
that if $D$ has at least five edges, then $C$ has at least five
edges. Other examples, according to what is shown in Section~2,
are obtained by restricting ourselves to big cocycles such that
both of the subgraphs $G'$ and $G''$ have at least two vertices,
or to big cocycles such that one of $G'$ and $G''$ has a subgraph
that is a cycle. These other restrictions, as the previous one,
accord with the Lemma of Section~3. The last example, involving
cycles, is related to a characterization of $K_{3,3}$ and $K_5$
that may be found in \cite{M97} (Lemma~3; see also~\cite{S06},
Lemma on the Kuratowski Graphs (2)).

These restricted variants of our criterion might be interesting
from an algorithmic point of view. They may shorten a procedure
for checking planarity.

\vspace{2ex}

\noindent {\small \emph{Acknowledgement.} We are grateful to
Sergei Melikhov for reading this note and making interesting
comments. We are also grateful to Rade \v Zivaljevi\' c for
leading us to \cite{M84}, and for giving us the opportunity to
discuss the matters of this note. Our work was supported by the
Ministry of Science of Serbia (Grant ON174026).}

\end{document}